\documentclass[a4paper,12pt]{article}
\usepackage{amsfonts, amsmath}
\makeatletter
\newbox\bk@bxb
\newbox\bk@bxa
\newif\if@bkcont
\newcount\bk@lcnt

\def\breakboxskip{2pt}
\def\breakboxparindent{1.8em}

\def\breakbox{\vskip\breakboxskip\relax
\setbox\bk@bxb\vbox\bgroup
\advance\linewidth -2\fboxrule
\hsize\linewidth\@parboxrestore
\parindent\breakboxparindent\relax}

\def\bk@split{%
\@tempdimb\ht\bk@bxb 
\advance\@tempdimb\dp\bk@bxb
\setbox\bk@bxa\vsplit\bk@bxb to\z@ 
\setbox\bk@bxa\vbox{\unvbox\bk@bxa}
\setbox\@tempboxa\vbox{\copy\bk@bxa\copy\bk@bxb}
\advance\@tempdimb-\ht\@tempboxa
\advance\@tempdimb-\dp\@tempboxa}

\def\bk@addfsepht{%
\setbox\bk@bxa\vbox{\vskip\fboxsep\box\bk@bxa}}

\def\bk@addskipht{%
\setbox\bk@bxa\vbox{\vskip\@tempdimb\box\bk@bxa}}

\def\bk@addfsepdp{%
\@tempdima\dp\bk@bxa
\advance\@tempdima\fboxsep
\dp\bk@bxa\@tempdima}

\def\bk@addskipdp{%
\@tempdima\dp\bk@bxa
\advance\@tempdima\@tempdimb
\dp\bk@bxa\@tempdima}

\def\bk@line{%
\hbox to \linewidth{%
\hskip-2\fboxsep\vrule \@width\fboxrule\hskip.5\fboxsep\vrule \@width\fboxrule\hskip1.5\fboxsep
\box\bk@bxa\hfil
}}%

\def\endbreakbox{\egroup
\ifhmode\par\fi{\noindent\bk@lcnt\@ne
\@bkconttrue\baselineskip\z@\lineskiplimit\z@
\lineskip\z@\vfuzz\maxdimen
\bk@split\bk@addfsepht\bk@addskipdp
\ifvoid\bk@bxb 
\def\bk@fstln{\bk@addfsepdp
\hskip-\parindent\vbox{\llap{\raisebox{-2ex}{\rule{1.5\fboxsep}{\fboxrule}\hskip.5\fboxsep}}\bk@line\llap{\rule{1.5\fboxsep}{\fboxrule}\hskip.5\fboxsep}}}

\else 
\def\bk@fstln{\vbox{\llap{\raisebox{-2ex}{\rule{1.5\fboxsep}{\fboxrule}\hskip.5\fboxsep}}\bk@line}\hfil%
\advance\bk@lcnt\@ne
\loop
\bk@split\bk@addskipdp\leavevmode
\ifvoid\bk@bxb 
\@bkcontfalse\bk@addfsepdp
\vtop{\bk@line\noindent\hskip-2\fboxsep{\rule{1.5\fboxsep}{\fboxrule}}}%

\else 
\bk@line
\fi
\hfil\advance\bk@lcnt\@ne
\if@bkcont\repeat}%
\fi
\leavevmode\bk@fstln\par}\vskip\breakboxskip\relax}

\def\smp{\smallskip\par}

\def\pf{\noindent{\bf Proof~:}\ }
\def\findemo{~\leaders\hbox to 1em{\hss\  \hss}\hfill~\raisebox{.5ex}{\framebox[1ex]{}}\smp}

\def\mpn{\medskip\par\noindent}
\def\smpn{\smallskip\par\noindent}

\def\smp{\smallskip\par}
\def\smpn{\smallskip\par\noindent}

\def\mpoint{\;\;.}
\def\mvirg{\;\;,}

\def\Ind{{\rm Ind}}

\def\Hom{{\rm Hom}}

\def\Aut{{\rm Aut}}

\def\Id{{\rm Id}}

\def\op{^{op}}

\newcommand{\romain}[1]{\uppercase\expandafter{\romannumeral #1}}

\newcommand{\gMod}[1]{#1{\hbox{-}\mathsf{Mod}}}

\def\op{^{op}}

\newenvironment{enonce}[1]{\pagebreak[2]\refstepcounter{subsection}\refstepcounter{prop}\smpn{{\bf \thesection.\arabic{prop}.\ \ #1~:}}\begin{it} }{\end{it}\smp}
\newenvironment{enonce*}[1]{\pagebreak[2]\smpn{#1~:}\begin{it} }{\end{it}\smp}
\newcommand{\result}[1]{\begin{enonce}{#1}}
\def\fresult{\end{enonce}}

\newenvironment{mth}[1]{\begin{breakbox}\begin{enonce}{#1}}{\end{enonce}\end{breakbox}}
\newenvironment{mth*}[1]{\begin{breakbox}\begin{enonce*}{#1}}{\end{enonce*}\end{breakbox}}
\newenvironment{rem}[1]{\refstepcounter{subsection}\refstepcounter{prop} \mpn{{\bf \thesection.\arabic{prop}.}\ \ \bf#1\ :}}{\smp}

\def\dom{\backslash}
\makeatletter

\renewenvironment{itemize}{\ifnum \@itemdepth >3 \@toodeep\else \advance\@itemdepth \@ne
\edef\@itemitem{labelitem\romannumeral\the\@itemdepth}%
\list{\csname\@itemitem\endcsname}{\setlength{\topsep}{1ex}\setlength{\itemsep}{0pt}\def\makelabel##1{\hss\llap{##1}}}\fi}
{\endlist}
\makeatother
\makeatletter
\def\@sect#1#2#3#4#5#6[#7]#8{\ifnum #2>\c@secnumdepth
    \let\@svsec\@empty\else
    \refstepcounter{#1}\edef\@svsec{\csname the#1\endcsname .\hskip .5em}\fi
    \@tempskipa #5\relax
     \ifdim \@tempskipa>\z@
       \begingroup #6\relax
         \@hangfrom{\hskip #3\relax\@svsec}{\interlinepenalty \@M #8\par}%
       \endgroup
      \csname #1mark\endcsname{#7}\addcontentsline
        {toc}{#1}{\ifnum #2>\c@secnumdepth \else
                     \protect\numberline{\csname the#1\endcsname}\fi
                   #7}\else
       \def\@svsechd{#6\hskip #3\relax  
                  \@svsec #8\csname #1mark\endcsname
                     {#7}\addcontentsline
                          {toc}{#1}{\ifnum #2>\c@secnumdepth \else
                            \protect\numberline{\csname the#1\endcsname}\fi
                      #7}}\fi
    \@xsect{#5}}
\def\section{\@startsection {section}{1}{\z@}{-3.5ex plus-1ex minus
    -.2ex}{2.3ex plus.2ex}{\reset@font\Large\bf}}  

\makeatother
\renewenvironment{equation}{\refstepcounter{subsection}\refstepcounter{prop}$$}{\leqno{\bf (\theprop)}$$}

\def\mar[#1]{\ar@{-}[#1]|-{\object@{<}}}
\def\marb[#1]{\ar@{-}[#1]|{\object+{  }}}

\newcommand{\categ}[1]{{\langle#1\rangle}}
\def\Sets{{\mathsf{Sets}}}
\def\Fun{{\mathsf{Fun}}}
\begin{document}
\centerline{\Large\bf Bisets as categories,}\vspace{.2cm}
\centerline{\Large\bf and tensor product of induced bimodules}
\vspace{.5cm}
\centerline{Serge Bouc}
\vspace{1cm}
{\footnotesize {\bf Abstract~:} Bisets can be considered as categories. This note uses this point of view to give a simple proof of a Mackey-like formula expressing the tensor product of two induced bimodules. \vspace{1ex}\par
{\bf AMS subject classification (2000)~:} 16D20, 20C20.\par
{\bf Keywords~:} Biset, induced, bimodule, tensor product.
}
\section{Introduction}
Let $R$ be a commutative ring, let $G$ and $H$ be finite groups, let $X$ be a subgroup of $H\times G$, and $M$ be an $RX$-module. If $m\in M$ and $(h,g)\in X$, set $h\cdot m\cdot g^{-1}=(h,g)\cdot m$~: this is a slight extension of the usual correspondence between $R(H\times G)$-modules and $(RH,RG)$-bimodules.\par
The object of this note is to give a simple proof of the following result~:
\begin{mth}{{\bf Theorem}} \label{formula}Let $R$ be a commutative ring, let $G$, $H$, and $K$ be finite groups, let $X$ be a subgroup of $H\times G$ and $Y$ be a subgroup of $K\times H$. Let $M$ be an $RX$-module, and $N$ be an $RY$-module. Then there is an isomorphism of $(RK,RG)$-bimodules
$$(\Ind_{Y}^{K\times H}N)\otimes_{RH}(\Ind_{X}^{H\times G}M)\cong\!\!\!\!\!\!\mathop{\oplus}_{t\in [p_2(Y)\dom H/p_1(X)]}\!\!\!\!\!\!\Ind_{Y*{^{(t,1)}}X}^{K\times G}(N\otimes_{k_2(Y)\cap{^{t}k_1(X)}}{^{(t,1)}M})\mvirg$$
where the notation is as follows (cf.~\cite{doublact})~:
$$p_1(X)=\{h\in H\mid \exists g\in G,\;(h,g)\in X\},\;\;k_1(X)=\{h\in H\mid (h,1)\in X\}$$
$$p_2(Y)=\{h\in H\mid \exists k\in K,\;(k,h)\in Y\},\;\;k_2(Y)=\{h\in H\mid (1,h)\in Y\}$$
$$Y*{^{(t,1)}X}=\{(k,g)\in K\times G\mid \exists h\in H,\;(k,h)\in Y,\;(h^t,g)\in X\}\mpoint$$
The action of $(k,g)\in Y*{^{(t,1)}X}$ on $N\otimes_{k_2(Y)\cap{^{t}k_1(X)}}{^{(t,1)}M}$ is given by
$$k\cdot(n\otimes m)\cdot g^{-1}=(k\cdot n\cdot h^{-1})\otimes (h^t\cdot m\cdot g^{-1})\mvirg$$
if $h\in H$ is chosen such that $(k,h)\in Y$ and $(h^t,g)\in X$.
\end{mth}
\section{Functors over bisets}
Recall that when $G$ and $H$ are groups, {\em an $(H,G)$-biset} $U$ is a set equipped with a left action of $H$ and a right action of $G$ which commute, i.e. such that $(hu)g=h(ug)$ for any $h\in H$, $u\in U$, and $g\in G$.
\begin{mth}{Notation} Let $G$ and $H$ be groups. When $U$ is an $(H,G)$-biset, let $\categ{U}$ denote the following category~: 
\begin{itemize}
\item The objects of $\categ{U}$ are the elements of $U$.
\item If $u,v\in U$, then
$$\Hom_\categ{U}(u,v)=\{(h,g)\in H\times G\mid hu=vg\}\mpoint$$
\item If $u,v,w\in U$, the composition of the morphisms $(h,g):u\to v$ and $(h',g'):v\to w$ is the morphism $(h'h,g'g):u\to w$.
\item If $u\in U$, the identity morphism $\Id_u:u\to u$ is the pair $(1,1)\in G\times G$.
\end{itemize}
\end{mth}
Note that the category $\categ{U}$ is a groupoid (any morphism is an isomorphism), and that for any $u\in U$, the group 
$$A(u)=\Hom_\categ{U}(u,u)=\{(h,g)\in H\times G\mid hu=ug\}$$ 
is a subgroup of $H\times G$.\par
A functor $M$ from $\categ{U}$ to a category $\mathcal{C}$ consists of a collection of objects $M(u)$ of $\mathcal{C}$, for $u\in U$, together with morphisms 
$$M(h,g):M(u)\to M(hug^{-1})$$
in the category $\mathcal{C}$, for $(h,g)\in H\times G$, fulfilling the usual functorial conditions. In particular, for each $u\in U$, there is a group homomorphism $A(u)\to \Aut_{\mathcal{C}}M(u)$. \par
Functors from $\categ{U}$ to $\mathcal{C}$ are the objects of a category $\Fun(\categ{U},\mathcal{C})$, in which the morphisms are natural transformation of functors. 
\mpn
\begin{mth}{Notation} When $\mathcal{C}$ is a subcategory of the category $\Sets$ of sets, and $M$ is a functor $\categ{U}\to \mathcal{C}$, the image of $m\in M(u)$ by the map $M(h,g): M(u)\to M(hug^{-1})$, for $(h,g)\in H\times G$, will be denoted by $hmg^{-1}$.
\end{mth}
In this case, a functor $M:\categ{U}\to \mathcal{C}$ is a collection of objects $M(u)$ of $\mathcal{C}$, for $u\in U$, together with morphisms $m\mapsto hmg^{-1}: M(u)\to M(hug^{-1})$ in $\mathcal{C}$, for $(h,g)\in H\times G$, such that $h'(hmg^{-1})g'^{-1}=(h'h)m(g'g)^{-1}$ and $1m1=m$, for any $(h,g)$, $(h',g')$ in $H\times G$, any $u\in U$, and any $m\in M(u)$. \par
\begin{rem}{Example} Suppose that $\mathcal{C}=\Sets$. Then the disjoint union $\bigsqcup M=\mathop{\bigsqcup}_{u\in U}\limits M(u)$ becomes an $(H,G)$-biset, and the map $\bigsqcup M\to U$ sending elements of $M(u)$ to $u$, for $u\in U$, is a map of $(H,G)$-bisets. Conversely, if $\pi:S\to U$ is a map of $(H,G)$-bisets, then the assignment $u\mapsto \pi^{-1}(u)$ is a functor from $\categ{U}$ to $\Sets$.\par
In other words, a functor $\categ{U}\to\Sets$ is just an $(H,G)$-biset over $U$. More precisely, the category $\Fun(\categ{U},\Sets)$ of such functors is equivalent to the category of $(H,G)$-bisets over $U$.
\end{rem}
\begin{rem}{Example} \label{exmp}Let $R$ be a commutative ring. In the remainder of this note, the category $\mathcal{C}$ will be the category $\gMod{R}$ of (left) $R$-modules. If $M$ is functor from $\categ{U}$ to $\gMod{R}$, then for each $u\in U$, the $R$-module $M(u)$ has a natural structure of $RA(u)$-module.\par
Conversely, let $[H\dom U/G]$ be a set of representatives of $(H,G)$-orbits on $U$. Equivalently $[H\dom U/G]$ is a set of representatives of isomorphism classes in the category $\categ{U}$. Since $\categ{U}$ is a groupoid, it is equivalent to its full subcategory $[H\dom U/G]$. In particular, this yields an equivalence of categories
\begin{equation}\label{eqv}
\Fun(\categ{U},\gMod{R})\cong \prod_{u\in[H\dom U/G]}\gMod{RA(u)}\mpoint
\end{equation}
\end{rem}
\begin{rem}{Remark}\label{rem1}
In the situation of Example~\ref{exmp}, the direct sum
$$\Sigma(M)=\mathop{\oplus}_{u\in U}M(u)$$
has a natural structure of $(RH,RG)$-bimodule, i.e. using the usual group isomorphism $(h,g)\mapsto (h,g^{-1})$ from $H\times G\op$ to $H\times G$, of left $R(H\times G)$-module. \par
Moreover, is is easy to see that there is an isomorphism of $(RH,RG)$-bimodules
$$\Sigma(M)\cong\mathop{\oplus}_{u\in [H\dom U/G]}\Ind_{A(u)}^{H\times G}M(u)\mpoint$$
\end{rem}
\section{Product of bisets, and product of functors}
Let $G$, $H$ and $K$ be groups. If $U$ is an $(H,G)$-biset and $V$ is a $(K,H)$-biset, recall that the product (or {\em composition}) of $V$ and $U$ is the set
$$V\times_HU=(V\times U)/H\mvirg$$
where the right action of $H$ on $(V\times U)$ is defined by $(v,u)\cdot h=(vh,h^{-1}u)$, for $v\in V$, $u\in U$, and $h\in H$. The set $V\times_HU$ is a $(K,G)$-biset for the following action
$$\forall z\in K,\;\forall x\in G,\;\forall v\in V,\;\forall u\in U,\;\;z\cdot(v,_{_H}u)\cdot x=(zv,_{_H}ux)\mvirg$$
where $(v,_{_H}u)$ denotes the $H$-orbit of $(v,u)$. 
\begin{mth}{Definition} Let $G$, $H$, and $K$ be finite groups. Let $U$ be a  finite $(H,G)$-biset, and $V$ be a finite $(K,H)$-biset. If $M$ is a functor $\categ{U}\to\gMod{R}$ and $N$ is a functor $\categ{V}\to \gMod{R}$, the tensor product $N\otimes_HM$ is the functor $\categ{V\times_H U}\to \gMod{R}$ defined by
$$(N\otimes_HM)(v,_{_H}u)=\big(\mathop{\oplus}_{h\in H}N(vh)\otimes_RM(h^{-1}u)\big)/\mathcal{I}_{v,u}\mvirg$$
where $\mathcal{I}_{v,u} $ is the $R$-submodule generated by the elements of the form $$[ny\otimes y^{-1}m]_{hy}-[n\otimes m]_h\mvirg$$
where $y\in H$, and where $[n\otimes m]_h$ denotes the element $n\otimes m$ of the component indexed by $h\in H$ in the direct sum, for $n\in N(vh)$, and $m\in M(h^{-1}u)$. \par
If $(k,g)\in K\times G$, then by definition
$$k \,[n\otimes m]_h \,g=[kn\otimes mg]_h\mpoint$$
\end{mth}
\begin{rem}{Remark} \label{rem2}It follows from this definition that
$$(N\otimes_HM)\big((v,_{_H}u)\big)\cong N(v)\otimes_{RH_{v,u}}M(u)\mvirg$$
where $H_{v,u}$ is the set of elements $h\in H$ such that $vh=v$ and $hu=u$.
\end{rem}
\begin{mth}{Lemma} There is an isomorphism of $(RK,RG)$-bimodules
$$\Sigma(N)\otimes_{RH}\Sigma(M)\cong\Sigma(N\otimes_HM)\mvirg$$
sending (from right to left) the element $[n\otimes m]_h$ to $n\otimes_{RH}m$.
\end{mth}
\pf To be more precise, the map $\alpha$ from
$$\Sigma(N\otimes_HM)=\mathop{\oplus}_{(v,_{_H}u)\in V\times_HU} \big(\oplus_{h\in H} N(vh)\otimes_RM(h^{-1}u)\big)/\mathcal{I}_{v,u}$$
sending the element $[n\otimes m]_h$ in the component indexed by $(v,_{_H}u)$ to the element $n\otimes m$ of the tensor product
$$\Sigma(N)\otimes_{RH}\Sigma(M)=\big(\mathop{\oplus}_{v\in V}N(v)\big)\otimes_{RH}\big(\mathop{\oplus}_{u\in U}M(u)\big)$$
is well defined. To show that it is an isomorphism, define a map 
$$\beta: \Sigma(N)\otimes_{RH}\Sigma(M)\to\Sigma(N\otimes_HM)$$
in the following way~: choose a set $S$ of representatives of the classes $(v,_{_H}u)$. Now map the element $n\otimes_{RH}m\in N(v)\otimes M(u)\subseteq \Sigma(N)\otimes_{RH}\Sigma(M)$ to $[n\otimes m]_h$, where $h\in H$ is chosen such that $(vh^{-1},hu)\in S$. Again, it is easy to see that this map is well defined, and that the maps $\alpha$ and $\beta$ are mutual inverse isomorphisms of $(RK,RG)$-bimodules.\findemo
\begin{mth}{Corollary} Let $G$, $H$, and $K$ be finite groups. Let $X$ be a subgroup of $H\times G$ and $Y$ be a subgroup of $K\times H$. Let $M$ be an $RX$-module, and $N$ be an $RY$-module. Then there is an isomorphism of $(RK,RG)$-bimodules
$$(\Ind_{Y}^{K\times H}N)\otimes_{RH}(\Ind_{X}^{H\times G}M)\cong\!\!\!\!\!\!\mathop{\oplus}_{t\in [p_2(Y)\dom H/p_1(X)]}\!\!\!\!\!\!\Ind_{Y*{^{(t,1)}}X}^{K\times G}(N\otimes_{k_2(Y)\cap{^{t}k_1(X)}}{^{(t,1)}M})\mpoint$$
\end{mth}
\pf Set $U=(H\times G)/X$. Then $U$ is an $(H,G)$-biset by $h\cdot(t,s)X\cdot g=(ht,g^{-1}s)X$, and this biset is transitive. If $u$ is the point $X$ of $U$, then $A(u)=X$, and the equivalence of categories~\ref{eqv} reads
$$\Fun(\categ{U},\gMod{R})\cong \gMod{RX}\mpoint$$
More precisely, for an $RX$-module $M$, this equivalence yields a functor $\tilde{M}:\categ{U}\to\gMod{R}$ in the following way~: for any $(h,g)\in H\times G$, set
$$\tilde{M}\big((h,g)X\big)=M\mpoint$$
Next, fix a set $S$ of representatives of elements of $U$, i.e. $X$-cosets in $H\times G$. For $(t,s)\in S$, and $(h,g)\in H\times G$, define a map
$$\tilde{M}(h,g): \tilde{M}\big((t,s)X\big)=M\to\tilde{M}\big((ht,gs)X\big)=M$$
by $\tilde{M}(h,g)(m)=(y,x)m$, where $(y,x)$ is the unique element of $X$ such that $(ht,gs)(y,x)^{-1}\in S$.\par
Then it is easy to check that $\tilde{M}$ is indeed a functor, and that there is an isomorphism of $(RH,RG)$-bimodules 
$$\Sigma(\tilde{M})\cong \Ind_X^{H\times G}M\mpoint$$
Similarly, set $V=(K\times H)/Y$, and define a functor $\tilde{N}\categ{V}\to\gMod{R}$, using the $RY$-module $N$. Then the corollary is a straightforward consequence of the lemma, applied to the functors $\tilde{M}$ and $\tilde{N}$, using Remark~\ref{rem1} and Remark~\ref{rem2}.\findemo

\begin{thebibliography}{1}

\bibitem{doublact}
S.~Bouc.
\newblock Foncteurs d'ensembles munis d'une double action.
\newblock {\em J. of Algebra}, 183(0238):664--736, 1996.

\end{thebibliography}

\vspace{1cm}
Serge Bouc\\
CNRS-LAMFA\\
Universit\'e de Picardie\\
33 rue St Leu\\
80039 - Amiens Cedex 1\\
FRANCE\\
{\tt email : serge.bouc@u-picardie.fr}
\end{document}